\documentclass[11pt]{article}
\usepackage{amssymb}
\usepackage{amscd}
\input amssym.def
\input amssym  

\newtheorem{theorem}{Theorem}[section]
\newtheorem{corollary}{Corollary}[section]
\newtheorem{lemma}{Lemma}[section]

\newtheorem{remark}{Remark}[section]

\usepackage{amssymb} 
\title{Principal
subspaces of higher-level standard
$\widehat{\mathfrak{sl}(3)}$-modules} \author{ Corina Calinescu
\footnote{The author gratefully acknowledges partial support from NSF
grant DMS 0100495 and from the Center for Discrete Mathematics and
Theoretical Computer Science (DIMACS), Rutgers University.} } \date{}
\begin{document}
\maketitle

\renewcommand{\theequation}{\thesection.\arabic{equation}}
\renewcommand{\thetheorem}{\thesection.\arabic{theorem}}
\setcounter{equation}{0}
\setcounter{theorem}{0}
\setcounter{section}{0}

\begin{abstract}
We use the theory of vertex operator algebras and intertwining
operators to obtain systems of $q$-difference equations satisfied by
the graded dimensions of the principal subspaces of certain level $k$
standard modules for $\widehat{\goth{sl}(3)}$. As a consequence we
establish new formulas for the graded dimensions of the principal
subspaces corresponding to the highest weights $i
\Lambda_1+(k-i)\Lambda_2$, where $1 \leq i \leq k$ and $\Lambda_1$ and
$\Lambda_2$ are fundamental weights of $\widehat{\goth{sl}(3)}$.
\end{abstract}

\section{Introduction}

The theory of vertex operator algebras (\cite{B}, \cite{FLM};
cf. \cite{LL}) leads to classical and new combinatorial identities
(see e.g. \cite{LW2}-\cite{LW4}, \cite{Ca}, \cite{MP1} and \cite{MP2})
and $q$-difference equations (recursions) satisfied by the graded
dimensions (characters) of certain substructures of the standard
representations of affine Lie algebras (see \cite{CLM1}-\cite{CLM2}
and \cite{C1}).

This paper is a continuation of \cite{C1}, to which we refer the
reader for background and notation (see also \cite{C2}). In \cite{C1}
we have derived a complete set of recursions that characterize the
graded dimensions of all the principal subspaces of the level $1$
standard representations of $\widehat{\goth{sl}(n)}$ with $n \geq
3$. Here we extend this approach to the principal subspaces of the
higher-level standard modules for $\widehat{\goth{sl}(3)}$. This work
and the work done in \cite{C1} can be viewed as a continuation of a
program to obtain Rogers-Ramanujan-type recursions, which was
initiated by Capparelli, Lepowsky and Milas in
\cite{CLM1}-\cite{CLM2}.

In the present paper we continue the study of a relationship between
intertwining operators, in the sense of \cite{FHL} and \cite{DL},
associated to standard modules and the corresponding principal
subspaces of these modules. We consider the definition of the
principal subspace $W(\Lambda)$ from \cite{FS1} and \cite{FS2},
namely, $W(\Lambda)=U(\bar{\goth{n}}) \cdot v_{\Lambda}$, where
$\Lambda$ and $v_{\Lambda}$ are the highest weight and a highest
weight vector of the standard $\widehat{\goth{sl}(3)}$-representation
$L(\Lambda)$ of level $k >1$. By $\bar{\goth{n}}$ we mean $\goth{n}
\otimes \mathbb{C}[t, t^{-1}]$, where $\goth{n}$ is the subalgebra of
$\goth{sl}(3)$ consisting of the strictly upper-triangular matrices.
As is well known, the highest weights of the level $k$ standard
modules have the form
$\Lambda=k_0\Lambda_0+k_1\Lambda_1+k_2\Lambda_2$, where $\Lambda_0,
\Lambda_1, \Lambda_2$ are the fundamental weights of
$\widehat{\goth{sl}(3)}$ and $k_0, k_1, k_2$ are nonnegative integers
whose sum is $k$.
 
The main result of our paper is the following theorem, which gives two
 families of exact sequences that yield linear systems of
 $q$-difference equations: \\ \\ {\bf Theorem } {\it For any integer
 $i$ with $1 \leq i \leq k$ there are natural sequences
\begin{eqnarray} \nonumber
\lefteqn{0 \longrightarrow W(i \Lambda_1 + (k-i) \Lambda_2)
\longrightarrow} \nonumber \\ && W(i \Lambda_0 + (k-i) \Lambda_1)
\longrightarrow \nonumber \\ && \hspace{2em} W((i-1)\Lambda_0+(k-i+1)
\Lambda_1) \longrightarrow 0 \nonumber
\end{eqnarray}
and
\begin{eqnarray} \nonumber
\lefteqn{0 \longrightarrow W((k-i) \Lambda_1 + i \Lambda_2)
\longrightarrow } \nonumber \\ && W(i \Lambda_0 + (k-i) \Lambda_2)
\longrightarrow \nonumber \\ && \hspace{2em} W((i-1)
\Lambda_0+(k-i+1)\Lambda_2) \longrightarrow 0, \nonumber
\end{eqnarray}
and these sequences are exact.}   
(See Theorem \ref{theorem1} for details.)
\vspace{1em}

As a consequence of this theorem we derive a system of $q$-difference
 equations that characterize the graded dimensions of the principal
 subspaces $W(i\Lambda_0+(k-i)\Lambda_j)$, where $0 \leq i \leq k$ and
 $j=1,2$ (see Theorem \ref{q-equations} below). The graded dimensions
 of these subspaces were previously obtained by Georgiev using a
 different method. Combining our recursions with the formulas for the
 graded dimensions of $W(i\Lambda_0+(k-i)\Lambda_j)$ we obtain the
 graded dimensions of the principal subspaces $W(i \Lambda_1+ (k-i)
 \Lambda_2)$ (see Corollary \ref{new} in Section 4). These are new
 results in the process of obtaining graded dimensions of principal
 subspaces; the method used in \cite{G1} did not give an answer for
 highest weights of type $i \Lambda_1+(k-i) \Lambda_2$ with $1 \leq i
 \leq k-1$.  Perhaps the main point of our strategy and results,
 though, is understanding the role of the intertwining operators
 involved in the construction of the exact sequences.
 
In this paper we formulate as a conjecture a presentation of the
principal subspaces of all the higher-level standard modules for
$\widehat{\goth{sl}(3)}$. It appears that in fact one can use an idea
developed in \cite{C1} to prove this result, but we will instead prove
this conjecture in a different way, in a future publication that is
part of ongoing joint work with Lepowsky and Milas.
   
The paper is organized as follows. Section 2 gives background and
notation. In Section 3 we discuss the principal subspaces of the
standard $\widehat{\goth{sl}(3)}$-modules. In Section 4 we construct
exact sequences that give $q$-difference equations and we obtain the
graded dimensions of the principal subspaces $W(i
\Lambda_1+(k-i)\Lambda_2)$ with $1 \leq i \leq k$.

This paper is part of the author's Ph.D. dissertation written under
the direction of James Lepowsky at Rutgers University.

{\bf Acknowledgement.} I would like to thank James Lepowsky for his
generous and inspiring guidance and for his encouragement during this
work. I am thankful to Antun Milas for many useful conversations.

\section{Preliminaries}

The aim of this section is to recall the vertex operator construction
of the higher-level standard modules for the untwisted affine Lie
algebra $\widehat{\goth{sl}(3)}$ by using the corresponding
contructions of the level $1$ standard representations. We use the
setting and notation from \cite{C1}.

In this paper we work with the finite-dimensional complex Lie algebra
$\goth{sl}(3)$ that has a standard basis
$$\{ x_{ \alpha _{1}}, \; x_{ \alpha _{2}}, \; x_{ \alpha _{1}+ \alpha
_{2}}, \; h_{\alpha_1}, \; h_{\alpha_2}, \; x_{- \alpha _{1}}, \; x_{-
\alpha _{2}}, \; x_{- \alpha_{1} - \alpha _{2}} \}.$$ We fix the
Cartan subalgebra $\goth{h}=\mathbb{C} h_{\alpha_1} \oplus
\mathbb{C}h_{\alpha_2}$ of ${\goth{sl}(3)}$. The standard symmetric
invariant nondegenerate bilinear form $\langle x,y
\rangle=\mbox{tr}(xy)$ defined for any $x$ and $y$ in ${\goth{sl}(3)}$
allows us to identify $\goth{h}$ with $\goth{h}^{*}$. Take $\alpha_1$
and $\alpha_2$ to be the (positive) simple roots corresponding to the
vectors $x_{\alpha_1}$ and $x_{\alpha_2}$. Under our identification we
have $h_{\alpha_1}=\alpha_1$ and $h_{\alpha_2}=\alpha_2$. The
fundamental weights of ${\goth{sl}(3)}$ are linear functionals
$\lambda_1$ and $\lambda_2$ in the dual space $\goth{h}^{*}$
($=\goth{h}$). They are determined by the conditions $\langle
\lambda_i, \alpha_j \rangle = \delta_{i,j}$ for $i, j =1,2$, so that
\[
\lambda_1= \frac{2}{3} \alpha_1+ \frac{1}{3} \alpha_2 \; \; \mbox{and}
\; \; \lambda_2= \frac{1}{3} \alpha_1+\frac{2}{3} \alpha_2.
\]

We denote by $\goth{n}$ the positive nilpotent subalgebra of
${\goth{sl}(3)}$,
$$ \goth{n}= \mathbb{C}x_{\alpha_1} \oplus \mathbb{C}x_{\alpha_2}
\oplus \mathbb{C}x_{\alpha_1+\alpha_2},
$$ 
which can be viewed as the subalgebra consisting of the strictly
upper-triangular matrices.

Now we consider the untwisted affine Lie algebra associated to
${\goth{sl}(3)}$,
\begin{equation}
\widehat{{\goth{sl}(3)}}= {\goth{sl}(3)} \otimes \mathbb{C}[t, t^{-1}]
\oplus \mathbb{C}c,
\end{equation}
where $c$ is a non-zero central element and 
\begin{equation}
[ x \otimes t^m, y \otimes t^n ] = [x, y] \otimes t^{m+n} + m\langle
x, y \rangle \delta _{m+n, 0} c
\end{equation}
for any $x, y \in {\goth{sl}(3)}$ and $m, n \in \mathbb{Z}$.  By
adjoining the degree operator $d$ ($[d, x \otimes t^m]=m$, $[d,c]=0$)
to the Lie algebra $\widehat{{\goth{sl}(3)}}$ one obtains the affine
Kac-Moody algebra $\widetilde{{\goth{sl}(3)}}=\widehat{{\goth{sl}(3)}}
\oplus \mathbb{C}d$ (cf. \cite{K}).  We introduce the following
subalgebras of $\widehat{{\goth{sl}(3)}}$:
\begin{equation} \label{n}
\bar{\goth{n}}= \goth{n} \otimes \mathbb{C}[t, t^{-1}] ,
\end{equation}
\begin{equation}
\bar{\goth{n}}^{+}=\goth{n} \otimes \mathbb{C}[t],
\end{equation}
\begin{equation}
\widehat{\goth{h}} = \goth{h} \otimes \mathbb{C}[t, t^{-1}] \oplus
\mathbb{C}c,
\end{equation}
and
\begin{equation}
\widehat {\goth{h}}_{\mathbb{Z}} = \coprod _{m \in \mathbb{ Z}
\setminus {0}} \goth{h} \otimes t^m \oplus \mathbb{C}c.
\end{equation}
The latter is a Heisenberg subalgebra of $ \widehat{{\goth{sl}(3)}}$
in the sense that its commutator subalgebra is equal to its center,
which is one-dimensional. The usual extension of the form $\langle
\cdot, \cdot \rangle $ to $\goth{h} \oplus \mathbb{C}c \oplus
\mathbb{C}d$ will be denoted by the same symbol ($\langle c, c
\rangle=0$, $\langle d, d \rangle =0$ and $\langle c, d \rangle =1$).
We will identify $\goth{h} \oplus \mathbb{C}c \oplus \mathbb{C}d$ with
its dual $(\goth{h} \oplus \mathbb{C}c \oplus \mathbb{C}d)^{*}$ via
this form.  The simple roots of $\widehat{{\goth{sl}(3)}}$ are $
\alpha _{0}$, $\alpha _{1}, \alpha _{2}$. We denote the fundamental
weights of $\widehat{{\goth{sl}(3)}}$ by $\Lambda_{0}, \Lambda _{1},
\Lambda _{2} $.  Then
$$
\alpha_0= c-(\alpha_1+\alpha_2)
$$
and
$$ \Lambda_0= d, \; \Lambda_1=\Lambda_0+\lambda_1, \;
\Lambda_2=\Lambda_0+\lambda_2.
$$
 
We say that a $\widehat{{\goth{sl}(3)}}$-module has level $k \in
 \mathbb{C}$ if $c$ acts as multiplication by $k$. It is well known
 that any standard module $L(\Lambda)$ with $\Lambda \in ( \goth{h}
 \oplus \mathbb{C}c \oplus \mathbb{C}d)^{*}$ has nonnegative integral
 level, given by $\langle \Lambda, c \rangle$ (cf. \cite{K}).  We
 denote by $L(\Lambda _{0})$, $ L( \Lambda _{1})$, $L(\Lambda_{2})$
 the standard $ \widehat {{\goth{sl}(3)}}$-modules of level 1 with
 $v_{\Lambda_0}$, $v_{\Lambda_1}$ and $v_{\Lambda_2}$ highest weight
 vectors.

We form the induced $ \widehat{\goth{h}}$-module
$$ M(1)= U( \widehat{\goth{h}}) \otimes _{U(\goth{h} \otimes
\mathbb{C}[t] \oplus \mathbb{C}c)}\mathbb{C},$$ such that $\goth{h}
\otimes \mathbb{C}[t]$ acts trivially and $c$ acts as identity on the
one-dimensional module $ \mathbb{C}$.  Let $Q= \mathbb{Z} \alpha _{1}
\oplus \mathbb{Z} \alpha _{2}$ be the root lattice and $P = \mathbb{Z}
\lambda _{1} \oplus \mathbb{Z} \lambda _{2}$ be the weight lattice of
${\goth{sl}(3)}$.  Denote by $\mathbb{C}[Q]$ and $\mathbb{C}[P]$ the
group algebras of the lattices $Q$ and $P$ with bases $\{ e^{\alpha} |
\alpha \in Q \}$ and $\{ e^{\lambda} | \lambda \in P \}$. Consider the
following vector spaces:
$$ V_P = M(1) \otimes \mathbb{C}[P],$$
$$V_Q= M(1) \otimes \mathbb{C}[Q]$$ 
and $$ V_Qe^{ \lambda_{i}}= M(1)
\otimes \mathbb{C}[Q]e^{ \lambda_{i}}, \; \; \; i=1,2.$$

It is well known (\cite{FK}, \cite{S}; cf. \cite{FLM}) that the vector
spaces $V_P$, $V_Q$, $V_Qe^{\lambda_1}$ and $V_Qe^{\lambda_2}$ admit a
natural $\widehat{{\goth{sl}(3)}}$-module structure via certain vertex
operators. Moreover, $V_Q$, $V_Qe^{\lambda_1}$ and $V_Qe^{\lambda_2}$
are the level $1$ standard representations of
$\widehat{{\goth{sl}(3)}}$ with highest weights $\Lambda_0$,
$\Lambda_1$ and $\Lambda_2$ and highest weight vectors
$v_{\Lambda_0}=1 \otimes 1$, $v_{\Lambda_1}= 1 \otimes e^{\lambda_1}$
and $v_{\Lambda_2}=1 \otimes e^{\lambda_2}$. We shall identify
\begin{equation}
V_Q \simeq L(\Lambda_0), \; V_Qe^{\lambda_1} \simeq L(\Lambda_1), \;
V_Qe^{\lambda_2} \simeq L(\Lambda_2)
\end{equation}
and we shall write 
\begin{equation} \label{vectors}
v_{\Lambda_0}=1, \; \; v_{\Lambda_1}=e^{\lambda_1} \; \; \mbox{and} \;
\; v_{\Lambda_2}=e^{\lambda_2}.
\end{equation}
For any weight $\lambda \in P$, by $e^{\lambda}$ we mean a vector of
 $V_P$ or an operator on $V_P$, depending on the context.  The space
 $V_Q$ has a natural structure of vertex operator algebra and
 $V_Qe^{\lambda_i}$ are $V_Q$-modules for $i=1,2$ (\cite{B} and
 \cite{FLM}).  See Section 2 in \cite{C1} for further details and
 background about the vertex operator construction of the level $1$
 standard $\widehat{{\goth{sl}(3)}}$-modules.
  
The highest weights of the level $k$ standard 
$\widehat{{\goth{sl}(3)}}$-modules are given by 
\begin{equation} \label{hw}
k_0 \Lambda_0 + k_1 \Lambda_1 + k_2 \Lambda_2,
\end{equation}
where $k_0, k_1, k_2$ are nonnegative integers such that
$k_0+k_1+k_2=k$. Throughout this paper $k>1$ stands for the level of a
representation.

We consider $L(\Lambda)$ a standard $\widehat{{\goth{sl}(3)}}$-module
of level $k$ with highest weight $\Lambda$ as in (\ref{hw}) and a
highest weight vector $v_{\Lambda}$.  Set
\begin{equation}
V_P^{\otimes k}= \underbrace{V_P \otimes \cdots \otimes V_P}_{k \; \;
\mbox{times}}
\end{equation}
and
\begin{equation}
v_{i_1, \dots, i_k}= v_{\Lambda_{i_1}} \otimes \cdots \otimes
v_{\Lambda_{i_k}} \in V_P^{\otimes k},
\end{equation}
where exactly $k_0$ indices are equal to $0$, $k_1$ indices are equal
to $1$ and $k_2$ indices are equal to $2$.  Then of course $v_{i_1,
\dots, i_k}$ is a highest weight vector for $\widehat{\goth{sl}(3)}$,
and
\begin{equation}
L(\Lambda) \simeq U(\widehat{{\goth{sl}(3)}}) \cdot v_{i_1, \dots,
i_k} \subset V_P ^{\otimes k}
\end{equation}
(cf. \cite{K}). 
Thus there is an embedding of $L(\Lambda)$ into $V_P^{\otimes k}$:
\begin{equation} \label{embedding}
L(\Lambda) \longrightarrow V_P^{\otimes k},
\end{equation}
uniquely determined by the identification $v_{\Lambda}=v_{i_1, \dots,
i_k}$.

The action of $\widehat{{\goth{sl}(3)}}$ on $V_P^{\otimes k}$ is given
by the usual comultiplication
\begin{equation} \label{action}
a \cdot v= \Delta(a) v= (a\otimes1 \otimes \cdots \otimes 1+ \cdots +
1 \otimes \cdots \otimes 1 \otimes a) v
\end{equation}
for $a \in \widehat{{\goth{sl}(3)}}$ and $v \in V_P^{\otimes k}$ and
this action extends in the usual way to $U(\widehat{{\goth{sl}(3)}})$.
Throughout this paper we will write $x(m)$ for the action of $x
\otimes t^m $ on any $\widehat{{\goth{sl}(3)}}$-module, where $x \in
{\goth{sl}(3)}$ and $m \in \mathbb{Z}$. In particular, for any root
$\alpha$ and integer $m$ we set $x_{\alpha}(m)$ for the action of
$x_{\alpha} \otimes t^m \in \widehat{{\goth{sl}(3)}}$ on
$L(\Lambda)$. In this paper we will also use the notation $x(m)$ for
the Lie algebra element $x \otimes t^m$. It will be clear from the
context if $x(m)$ is an operator or an element of
$\widehat{{\goth{sl}(3)}}$.

The standard $\widehat{{\goth{sl}(3)}}$-modules have structures of
vertex operator algebra and modules and this result is stated below:

\begin{theorem}(\cite{FZ}; cf. \cite{DL}, \cite{Li1} and \cite{LL}) 
The standard module $L(k \Lambda_0)$ has a natural vertex operator
algebra structure. The set of the level $k$ standard
$\widehat{{\goth{sl}(3)}}$-modules provides a complete list of
irreducible $L(k \Lambda_0)$-modules (up to equivalence).
\end{theorem}

The vertex operator map
\begin{eqnarray} 
Y(\cdot, x): L(k\Lambda_0) & \longrightarrow & \mbox{End} \; L(k\Lambda_0)
[[x, x^{-1}]] \nonumber \\
v & \mapsto & Y(v, x)\sum_{m \in \mathbb{Z}} v_m x^{-m-1}
\nonumber
\end{eqnarray}
has the property
\begin{equation} 
Y(x_{\alpha}(-1) \cdot v_{k\Lambda_0}, x)= \sum_{m \in \mathbb{Z}}
x_{\alpha}(m)x^{-m-1}
\end{equation}
(cf. \cite{FZ}, \cite{DL}, \cite{Li1} and \cite{LL}).

We have the following products of operators:
\begin{equation} \label{product oper}
e^{\lambda} x_{\alpha}(m)= x_{\alpha}(m- \langle \lambda, \alpha
\rangle ) e^{\lambda}
\end{equation}
for any $\lambda \in P$, $\alpha$ a simple root and $m \in \mathbb{Z}$
and
\begin{equation} \label{prod-opera}
e^{\lambda} h(n)= h(n) e^{\lambda}
\end{equation}
for any $\lambda \in P$, $h \in \goth{h}$ and $n \neq 0$.

\section{Principal subspaces}
\setcounter{equation}{0}

In \cite{FS1}-\cite{FS2} B. Feigin and A. Stoyanovsky introduced
certain substructures, called principal subspaces, of the standard
modules for $\widehat{\goth{sl}(n)}$ with $n \geq 2$. Consider
$L(\Lambda)$ a level $k$ standard module for
$\widehat{{\goth{sl}(3)}}$. Then $\Lambda=k_0\Lambda_0+k_1
\Lambda_1+k_2 \Lambda_2$, where $k_0, k_1, k_2 \in \mathbb{N}$ whose
sum is $k$ (cf. (\ref{hw})). The principal subspace of $L(\Lambda)$,
denoted by $W(\Lambda)$, is defined as follows:
\begin{equation}
W(\Lambda)= U(\bar{\goth{n}}) \cdot v_{\Lambda},
\end{equation}
where $U(\bar{\goth{n}})$ is the universal enveloping algebra of
$\bar{\goth{n}}$ (recall (\ref{n})) and $v_{\Lambda}$ is a highest
weight vector of $L(\Lambda)$.

Let us consider the natural surjective map:
\begin{eqnarray} \label{map}
f_{\Lambda} : U(\bar{\goth{n}}) & \longrightarrow & W(\Lambda)  \\
a & \mapsto & a \cdot v_{\Lambda}. \nonumber
\end{eqnarray}
Denote by $I_{\Lambda}$ the annihilator of the highest weight vector
$v_{\Lambda}$ in $U(\bar{\goth{n}})$,
\begin{equation} \label{I}
I_{\Lambda}= \mbox{Ker} \; f_{\Lambda}.
\end{equation}
This is a left ideal of $U(\bar{\goth{n}})$. Hence the principal
subspace $W(\Lambda)$ is isomorphic (as linear spaces) with the
quotient space $U(\bar{\goth{n}}) / I_{\Lambda}$.

Recall from Section 2 the operators $x_{\alpha}(m)$ for any root
$\alpha$ and $m \in \mathbb{Z}$.  We consider the following formal
infinite sums:
\begin{equation}
R_t^{[1]}=\sum_{m_1 + \cdots + m_{k+1}=t} x_{\alpha_1}(m_1) \cdots
x_{\alpha_1}(m_{k+1})
\end{equation}
and
\begin{equation}
R_t^{[2]}=\sum_{m_1+ \cdots + m_{k+1}=t}x_{\alpha_2}(m_1) \cdots
x_{\alpha_2}(m_{k+1})
\end{equation}
for any $t \in \mathbb{Z}$.  It will be convenient to truncate
$R_t^{[1]}$ and $R_t^{[2]}$ as follows:
\begin{equation} \label{R_t^j}
R_{t;m}^{[j]}= \sum_{\begin{array} {c} m_1, \dots, m_{k+1} \leq
m, \\ m_1+ \cdots + m_{k+1}=t \end{array}}x_{\alpha_j}(m_1) \cdots
x_{\alpha_j}(m_{k+1}),
\end{equation}
where $j=1,2$ and $m$ is a fixed (not necessarily negative) integer.
We shall often view $R_{t;m}^{[j]}$ as elements of the algebra
$U(\bar{\goth{n}})$, rather than as endomorphisms of a
$\widehat{{\goth{sl}(3)}}$-module. It will be clear from the context
when expressions such as (\ref{R_t^j}) are understood as elements of
$U(\bar{\goth{n}})$ or as operators. Denote by $\widetilde{U(\bar{\goth{n}})}$ the completion of $U(\bar{\goth{n}})$ in the sense of \cite{LW3} or \cite{MP1}. Let us denote by ${\mathcal J}$ the
two-sided ideal of $\widetilde{U(\bar{\goth{n}})}$ generated by $R_{t}^{[1]}$
and $R_{t}^{[2]}$ for $t \in \mathbb{Z}$.

In order to obtain exact sequences and $q$-difference equations in the
next section, it is important to have a description of the ideals
$I_{k_0 \Lambda_0+k_1\Lambda_1+k_2 \Lambda_2}$. This problem is
equivalent to finding a presentation of the principal subspaces $W(k_0
\Lambda_0+k_1 \Lambda_1+k_2 \Lambda_2)$.

The presentation of the principal subspaces of the standard modules
for $\widehat{\goth{sl}(2)}$ was proved in \cite{FS1}-\cite{FS2}. In
\cite{FS2} the authors announced a presentation of the principal
subspaces $W(k \Lambda_0)$ of the standard representations $L(k
\Lambda_0)$ of $\widehat{\goth{sl}(n)}$ for $n \geq 3$. We have given
a precise description of the ideals $I_{\Lambda}$ corresponding to the
level $1$ standard modules $L(\Lambda)$ of $\widehat{\goth{sl}(n)}$
for $n \geq 3$ in \cite{C1}.

Here we conjecture a description of the ideals
$I_{k_0\Lambda_0+k_1\Lambda_1+k_2\Lambda_2}$ as follows: \\
 
{\bf Conjecture} \label{presentation}
We have
\begin{equation}
I_{k\Lambda_0} \equiv {\mathcal J} \; \; \; \mbox{modulo} \; \; \; \widetilde{U(\bar{\goth{n}}) \bar{\goth{n}}^{+}}
\end{equation}
and
\begin{equation}
I_{k_0\Lambda_0+ (k-k_0)\Lambda_j} \equiv {\mathcal J} + U(\bar{\goth{n}})
x_{\alpha_j}(-1)^{k_0+1} \; \; \; \mbox{modulo} \; \; \; \widetilde{U(\bar{\goth{n}}) \bar{\goth{n}}^{+}}
\end{equation}
for $0 \leq k_0 \leq k$ and $j=1,2$. More generally, 
\begin{equation}
I_{k_0\Lambda_0+k_1\Lambda_1+k_2 \Lambda_2} \equiv  {\mathcal J} 
+ U(\bar{\goth{n}}) x_{\alpha_1}(-1)^{k_0+k_2+1} +
U(\bar{\goth{n}}) x_{\alpha_2}(-1)^{k_0+k_1+1} \; \; \;
\mbox{modulo} \; \; \;  \widetilde{U(\bar{\goth{n}}) \bar{\goth{n}}^{+}}
\end{equation}
for any $k_0, k_1, k_2 \geq 0$ such that $k_0+k_1+k_2=k$.
\\ \\

As a consequence of this statement we obtain the discrepancy between
the defining ideals of the principal subspaces.

\begin{corollary} \label{discrepancy}
We have
\begin{equation}
I_{k_0 \Lambda_0+(k-k_0) \Lambda_1}= I_{k \Lambda_0}
+U(\bar{\goth{n}}) x_{\alpha_1}(-1)^{k_0+1},
\end{equation}
\begin{equation}
I_{k_0 \Lambda_0+(k-k_0) \Lambda_2}= I_{k \Lambda_0}
+U(\bar{\goth{n}}) x_{\alpha_2}(-1)^{k_0+1}
\end{equation}
and
\begin{equation}
I_{k_0\Lambda_0+k_1\Lambda_1+k_2 \Lambda_2}= I_{k \Lambda_0}+
U(\bar{\goth{n}}) x_{\alpha_1}(-1)^{k_0+k_2+1} + U(\bar{\goth{n}})
x_{\alpha_2}(-1)^{k_0+k_1+1}.
\end{equation}
\end{corollary}

In order to prove the above conjecture we can follow the idea of the
proof of a presentation of the principal subspaces of the level $1$
standard modules for $\widehat{\goth{sl}(n)}$ for $n \geq 3$
developed in \cite{C1}. Since this proof of a presentation of the
principal subspaces of higher-level representations is very technical
we omit it. We will instead prove this conjecture in \cite{CalLM4} by
using a different and new approach.

\section{Exact sequences and $q$-difference equations}
\setcounter{equation}{0}

In this section we prove our main results. Mainly, we construct exact
sequences of maps between principal subspaces of certain standard
$\widehat{\goth{sl}(3)}$-modules. We derive systems of $q$-difference
equations satisfied by the graded dimensions of the principal
subspaces $W(i\Lambda_0+(k-i)\Lambda_j)$ for $0 \leq i \leq k$ and
$j=1,2$. We also obtain new formulas for the graded dimensions of
$W(i\Lambda_1+(k-i)\Lambda_2)$.

Take $\Lambda= k_0 \Lambda_0+k_1 \Lambda_1+ k_2 \Lambda_2$, where
$k_0, k_1, k_2$ are nonnegative integers with $k_0+k_1+k_2=k$ and take
$L(\Lambda)$ a standard module of level $k$.  It is known
(cf. \cite{FZ}, \cite {DL}, \cite{LL}) that $L(\Lambda)$ is graded
with respect to a standard action of the Virasoro algebra operator
$L(0)$ and that
\begin{equation}
L(\Lambda)= \coprod_{s \in \mathbb{Z}} L(\Lambda)_{s+ h_{\Lambda}},
\end{equation}
where
\begin{equation} \label{h}
h_{\Lambda}= \frac{\langle \Lambda, \Lambda+\alpha_1+\alpha_2
\rangle}{2(k+3)}
\end{equation}
and where $L(\Lambda)_{s+h_{\Lambda}}$ is the weight space of
$L(\Lambda)$ of weight $s+h_{\Lambda}$.  This is known as grading by
{\it weight}. The space $L(\Lambda)$ has also gradings by {\it charge}
given by the eigenvalues of the operators $\lambda_1$ and $\lambda_2$
(thought as $\lambda_1(0)$ and $\lambda_2(0)$). The gradings by charge
are compatible with the weight grading. Thus $L(\Lambda)$ decomposes
as
\begin{equation}
L(\Lambda)= \coprod_{r_1, r_2, s \in \mathbb{Z}} L(\Lambda)_{r_1+
\langle \lambda_1, \Lambda \rangle, r_2+ \langle \lambda_2, \Lambda
\rangle; s+ h_{\Lambda}},
\end{equation}
where $L(\Lambda)_{r_1+ \langle \lambda_1, \Lambda \rangle, r_2+
\langle \lambda_2, \Lambda \rangle; s+ h_{\Lambda}}$ is the subspace
of $L(\Lambda)$ consisting of the vectors of charges $r_1+\langle
\lambda_1, \Lambda \rangle$, $r_2+\langle \lambda_2, \Lambda \rangle$
and of weight $s+ h_{\Lambda}$.

Now we restrict these gradings to the principal subspace $W(\Lambda)
\subset L(\Lambda)$ and thus we have
\begin{equation}
W(\Lambda)= \coprod_{r_1, r_2, s \in \mathbb{N}} W(\Lambda)_{r_1+
\langle \lambda_1, \Lambda \rangle, r_2+ \langle \lambda_2, \Lambda
\rangle; s+ h_{\Lambda}}.
\end{equation} 
We consider the graded dimension (i.e. the generating function of the
dimensions of the homogeneous subspaces) of the space $W(\Lambda)$:
\begin{equation}
\chi_{W(\Lambda)}(x_1,x_2;q)=\mbox{dim}_{*} (W(\Lambda), x_1, x_2;
q)=\mbox{tr}|_{W(\Lambda)}x_1^{\lambda_1}x_2^{\lambda_2} q^{L(0)},
\end{equation}
where $x_1, x_2$ and $q$ are formal variables. To avoid the factor
$x_1^{\langle \lambda_1, \Lambda \rangle} x_2^{\langle \lambda_2,
\Lambda \rangle}q^{h_{\Lambda}}$ we use slightly modified graded
dimensions as follows:
\begin{equation}
\chi'_{W(\Lambda)}(x_1, x_2; q)= x_1^{-\langle \lambda_1, \Lambda
\rangle} x_2^{-\langle \lambda_2, \Lambda \rangle} q^{-h_{\Lambda}}
\chi_{W(\Lambda)} (x_1, x_2; q).
\end{equation} 
Thus we have
\[
\chi'_{W(\Lambda)}(x_1, x_2; q) \in \mathbb{C}[[x_1, x_2; q]],
\]
and in fact, the constant term of $\chi'_{W(\Lambda)}(x_1, x_2; q)$ is
$1$. Notice that
\begin{equation} \label{chi0}
\chi'_{W(k \Lambda_0)}(x_1, x_2; q)= \chi_{W(k\Lambda_0)}(x_1, x_2; q)
\in \mathbb{C}[[x_1, x_2; q]].
\end{equation}
We shall also use the notation
\begin{equation} \label{w'}
W(\Lambda)'_{r_1, r_2; k}= W(\Lambda)_{r_1+\langle \lambda_1, \Lambda
\rangle, r_2+\langle \lambda_2, \Lambda \rangle; k+h_{\Lambda}}.
\end{equation}

We now recall from \cite{C1} the automorphisms $e^{\lambda}: V_P
\longrightarrow V_P$ for any weight $\lambda$. Consider
\begin{equation} \label{e_k}
e^{\lambda}_{(k)} : V_P^{\otimes k} \longrightarrow V_P^{\otimes k},
\end{equation}
$$ 
e^{\lambda}_{(k)}=\underbrace{e^{\lambda} \otimes \cdots \otimes
e^{\lambda}}_{k \; \; \mbox{times}},
$$
a linear automorphism for any $\lambda \in P$. 
Thus it follows that
\begin{equation} \label{form1} 
e^{\lambda}_{(k)} \; (x_{\alpha}(m_{1}) \cdots x_{\alpha}(m_{r}) \cdot
v)= x_{\alpha}(m_1-\langle \lambda, \alpha \rangle) \cdots
x_{\alpha}(m_r- \langle \lambda, \alpha \rangle) \cdot
e^{\lambda}_{(k)} (v)
\end{equation}
and
\begin{equation} \label{form2}
e^{\lambda}_{(k)} \; (x_{-\alpha}(m_{1}) \cdots x_{-\alpha}(m_{r})
\cdot v)= x_{-\alpha}(m_1+\langle \lambda, \alpha \rangle) \cdots
x_{-\alpha}(m_r+ \langle \lambda, \alpha \rangle) \cdot
e^{\lambda}_{(k)} (v)
\end{equation}
for any positive root $\alpha$, $m_1, \dots , m_r \in \mathbb{Z}$ and
$v \in V_P^{\otimes k}$ (cf. (\ref{product oper})). We also have
\begin{equation} \label{form3}
e^{\lambda}_{(k)} \; (h_1(m_{1}) \cdots h_1(m_{r}) \cdot v)= h_1(m_1)
\cdots h_r(m_r) \cdot e^{\lambda}_{(k)}( v)
\end{equation}
if each $m_j \neq 0$ and $h_j \in \goth{h}$ (recall (\ref{prod-opera})).

One can also consider maps of the following type:
\begin{equation} \label{id}
{\rm Id}_{(k_0)} \otimes e^{\lambda}_{(k-k_0)}: V_P^{\otimes k}
\longrightarrow V_P^{\otimes k}
\end{equation}
for any $\lambda \in P$ and $0 \leq k_0 \leq k$, where $\mbox{Id}$ is
the identity map.

Let $k_0$ and $k_1$ be nonnegative integers whose sum is $k$. The next
results show that there exist maps of type (\ref{e_k}) between certain
standard modules, and in particular between their corresponding
principal subspaces.

\begin{lemma} \label{lemma1} 
Let $i$ be an integer such that $0 \leq i \leq k$.
\begin{enumerate}
\item The restriction of $e^{ \lambda_1}_{(k)}$ to
$L(i\Lambda_0+(k-i)\Lambda_2) $ lies in $L((k-i)\Lambda_0+i
\Lambda_1)$.
\item We have an injective linear map between principal subspaces 
\begin{equation} \label {principal}
e^{\lambda_1}_{(k)}: W(i \Lambda_0+(k-i)\Lambda_2) \longrightarrow
W((k-i)\Lambda_0 +i \Lambda_1).
\end{equation}
\item If $i=k$ then (\ref {principal}) is a linear isomorphism
\begin{equation} \label{iso1-0}
e^{\lambda_1}_{(k)}: W(k\Lambda_0) \longrightarrow W(k\Lambda_1).
\end{equation}
Moreover, we obtain the following relation between the graded
dimensions of the principal subspaces $W(k\Lambda_0)$ and
$W(k\Lambda_1)$:
\begin{equation} \label{1-0}
\chi'_{W(k \Lambda_1)}(x_1, x_2; q) = \chi'_{W(k \Lambda_0)} (x_1q,
x_2;q),
\end{equation}
which is equivalent with
\begin{equation}
\chi_{W(k \Lambda_1)}(x_1, x_2; q)
=x_1^{2k/3}x_2^{k/3}q^{h_{k\Lambda_1}} \chi _{W(k \Lambda_0)} (x_1q,
x_2; q).
\end{equation}
\end{enumerate}
\end{lemma}
{\em Proof:}
\begin{enumerate}
\item  
We view $L(i\Lambda_0+(k-i) \Lambda_2)$ and $L((k-i) \Lambda_0+i
\Lambda_1)$ embedded in $V_P^{\otimes k}$ (cf. (\ref{embedding})). We
have
$$L(i\Lambda_0 +(k-i)\Lambda_2) = U( \widehat{{\goth{sl}(3)}}) \cdot
v_{i_1, \dots, i_k},
$$
where $v_{i_1, \dots, i_k}=v_{\Lambda_{i_1}} \otimes \cdots \otimes
v_{\Lambda_{i_k}}$ with exactly $i$ indices equal to 0 and $k-i$
indices equal to 2. We may and do assume that the first $i$ indices
are 0 and the other $k-i$ indices are 2. By using our
identifications (\ref{vectors}) we have
$$
e^{\lambda_1}v_{\Lambda_2}=e^{\alpha_1+\alpha_2}=x_{\alpha_1+\alpha_2}
(-1) \cdot v_{\Lambda_0},
$$
$$ x_{\alpha_1+\alpha_2}(-1)\cdot v_{\Lambda_1} =
x_{\alpha_1+\alpha_2}(-1) \cdot e^{\lambda_1}=0
$$
and
$$ x_{\alpha_1+\alpha_2} (-1)^2 \cdot v_{\Lambda_0}
=x_{\alpha_1+\alpha_2}(-1)^2 \cdot 1=0
$$
(cf. \cite{G1} and \cite{C1}).
Thus it follows that
\begin{eqnarray} \label{form4}
&&e^{\lambda_1}_{(k)} (v_{i \Lambda_0+(k-i)
\Lambda_2})=e^{\lambda_1}_{(k)} (\underbrace{v_{\Lambda_0} \otimes
\cdots \otimes v_{\Lambda_0}}_{i\; \; \mbox{times}} \otimes
\underbrace{v_{\Lambda_2} \otimes \cdots \otimes
v_{\Lambda_2}}_{(k-i)\; \; \mbox{times}}) \\ \nonumber &&
=\underbrace{ v_{\Lambda_1} \otimes \cdots \otimes v_{\Lambda_1}}_{i
\; \; \mbox{times}} \otimes \underbrace{x_{\alpha_1 + \alpha_2} (-1)
\cdot v_{\Lambda_0} \otimes \cdots \otimes x_{\alpha_1+\alpha_2} (-1)
\cdot v_{\Lambda_0}}_{(k-i) \; \; \mbox{times}} \nonumber \\ &&=
a(\Delta (x_{\alpha_1+\alpha_2}(-1)^{k-i}) (\underbrace {v_{\Lambda_1}
\otimes \cdots \otimes v_{\Lambda_1}}_{i \; \; \mbox{times}} \otimes
\underbrace{v_{\Lambda_0} \otimes \cdots \otimes v_{\Lambda_0}}_{(k-i)
\; \; \mbox{times}}), \nonumber\\ && = a
x_{\alpha_1+\alpha_2}(-1)^{k-i} \cdot v_{(k-i)\Lambda_0+i
\Lambda_1},\nonumber
\end{eqnarray}
where $a$ is a nonzero constant, so that
\begin{equation} \label{form5}
e^{\lambda_1}_{(k)} ( v_{i \Lambda_0+(k-i) \Lambda_2}) \in
U(\widehat{{\goth{sl}(3)}}) \cdot v_{(k-i)\Lambda_0+i \Lambda_1}.
\end{equation}
Now by (\ref{form1}), (\ref{form2}), (\ref{form3}) and (\ref{form5})
we obtain the linear map
\begin{equation} \label{com1}
e^{\lambda_1}_{(k)} : L(i \Lambda_0 + (k-i) \Lambda_2) \longrightarrow
L((k-i) \Lambda_0 + i \Lambda_1).
\end{equation}
\item
Let us restrict the map (\ref{com1}) to the principal subspace
$W(i\Lambda_0+(k-i)\Lambda_2)$.  Using similar arguments with
$U(\bar{\goth{n}})$ instead of $U(\widehat{{\goth{sl}(3)}})$, we
obtain a linear map
\[ 
e^{\lambda_1}_{(k)} : W(i\Lambda_0 + (k-i) \Lambda_2) \longrightarrow
W((k-i)\Lambda_0 + i \Lambda_1),
\]
which is clearly injective.
\item
Let us take $i=k$ in (\ref{principal}). By (\ref{form4}) we have
\begin{equation}
e^{ \lambda_1}_{(k)} (v_{k \Lambda_0})=v_{k \Lambda_1}.
\end{equation} 
Since
$$ W(k \Lambda_1)=U(\bar{\goth{n}}) \cdot v_{k\Lambda_1}=
U(\bar{\goth{n}}) \cdot e^{\lambda_1}_{(k)} (v_{k\Lambda_0})=
e^{\lambda_1}_{(k)} (W(k \Lambda_0))
$$
(cf. (\ref{form1})) we obtain that the linear map
\begin{equation} \label{iso}
e^{\lambda_1}_{(k)}: W(k \Lambda_0) \longrightarrow W(k \Lambda_1)
\end{equation}
is surjective and thus it is a linear isomorphism. The isomorphism
 (\ref{iso}) does not preserve weight and charge.  Let $W(k
 \Lambda_0)_{r_1, r_2; s}$ with $r_1, r_2, s \in \mathbb{N}$ be an
 homogeneous subspace of $W(k \Lambda_0)$. The map (\ref{iso})
 increases the charge corresponding to $\lambda_j$ by $\langle k
 \Lambda_1, \lambda_j \rangle$ for $j=1,2$. For any $w \in W(k
 \Lambda_0)_{r_1, r_2; s}$ the homogeneous element
 $e^{\lambda_1}_{(k)} (w)$ has weight $s+r_1+h_{k \Lambda_1}$. Hence
 we obtain an isomorphism between homogeneous spaces
$$
e^{\lambda_1}_{(k)} : W(k\Lambda_0)'_{r_1, r_2; s} \longrightarrow 
W(k\Lambda_1)'_{r_1, r_2; s+r_1},
$$ 
which gives the relation between the graded dimensions 
$$
\chi'_{W(k\Lambda_1)}(x_1,x_2;q)=\chi'_{W(k\Lambda_0)}(x_1q,x_2;q),
$$
and so
$$ \chi_{W(k\Lambda_1)}(x_1, x_2; q) =
x_1^{2k/3}x_2^{k/3}q^{h_{k\Lambda_1}} \chi _{W(k\Lambda_0)} (x_1q,
x_2; q). \; \; \; \; \; \Box $$
\end{enumerate}
\vspace{1em}

We have a result completely analogous to the previous lemma:

\begin{lemma} \label{lemma2} 
Let $i$ be an integer such that $0 \leq i \leq k$.
\begin{enumerate}
\item The restriction of $e^{ \lambda_2}_{(k)}$ to $L(i
\Lambda_0+(k-i) \Lambda_1) $ lies in $L((k-i) \Lambda_0+i \Lambda_2)$.
\item At the level of the principal subspaces we have the following
injection:
\begin{equation} \label {bprincipal}
e^{\lambda_2}_{(k)}: W(i \Lambda_0+(k-i)\Lambda_1) \longrightarrow
W((k-i)\Lambda_0 +i \Lambda_2).
\end{equation}
\item If $i=k$ then (\ref {bprincipal}) is a linear isomorphism
\begin{equation} \label{iso2-0}
e^{\lambda_2}_{(k)}: W(k \Lambda_0) \longrightarrow W(k\Lambda_2).
\end{equation}
In particular, we obtain 
\begin{equation} \label{2-0}
\chi'_{W(k \Lambda_2)}(x_1, x_2; q) = \chi'_{W(k \Lambda_0)}(x_1,
x_2q; q),
\end{equation}
which is equivalent with
\begin{equation}
\chi_{W(k \Lambda_2)}(x_1, x_2; q) =
x_1^{k/3}x_2^{2k/3}q^{h_{k\Lambda_2}} \chi _{W(k\Lambda_0)} (x_1,
x_2q; q). \; \; \Box
\end{equation}
\end{enumerate}
\end{lemma}

Our main goal is to obtain exact sequences of maps between principal
subspaces. We have already seen in Lemmas \ref{lemma1} and
\ref{lemma2} that the maps $e^{\lambda_1}_{(k)}$ and
$e^{\lambda_2}_{(k)}$ are examples in this direction. In order to
construct more general exact sequences and thereby to obtain
$q$-difference equations we introduce maps of type $e^{\lambda}_{(k)}$
with $\lambda \neq \lambda_1$ and $\lambda \neq \lambda_2$ between
principal subspaces.

We first consider the weight $\lambda= \alpha_1-\lambda_1= \frac{1}{3}
\alpha_1- \frac{1}{3} \alpha_2$ and the linear isomorphism
\begin{equation} \label{1-2}
e^{\lambda}_{(k)}: V_P^{\otimes k} \longrightarrow V_P^{\otimes k}.
\end{equation}
The restriction of (\ref{1-2}) to the principal subspace $W(i
\Lambda_1+(k-i) \Lambda_2)$ is the map
\begin{equation} \label{map1}
e^{\lambda}_{(k)}: W(i \Lambda_1+(k-i) \Lambda_2) \longrightarrow W(i
\Lambda_0+ (k-i) \Lambda_1),
\end{equation}
where $0 \leq i \leq k$.

Since
$$
\lambda+\lambda_1=\alpha_1, \; \; \lambda+\lambda_2=\lambda_1
$$
and
$$ v_{i \Lambda_1+(k-i) \Lambda_2}=\underbrace{v_{\Lambda_1} \otimes
\cdots \otimes v_{\Lambda_1}}_{i \; \; \mbox{times}} \otimes
\underbrace{v_{\Lambda_2} \otimes \cdots \otimes v_{\Lambda_2}}_{(k-i)
\; \; \mbox{times}},
$$ 
\begin{eqnarray} \label{hwc}
&&e^{\lambda}_{(k)} (v_{i \Lambda_1+(k-i) \Lambda_2}) \\
&&=\underbrace{x_{\alpha_1}(-1) \cdot v_{\Lambda_0} \otimes \cdots
\otimes x_{\alpha_1}(-1) \cdot v_{\Lambda_0}}_{i \; \; \mbox{times}}
\otimes \underbrace{ v_{\Lambda_1} \otimes \cdots \otimes
v_{\Lambda_1}}_{(k-i) \; \; \mbox{times}} \nonumber \\ &&= a\Delta
(x_{\alpha_1}(-1)^i ) ( \underbrace{v_{\Lambda_0} \otimes \cdots
\otimes v_{\Lambda_0}}_{i \; \; \mbox{times}} \otimes \underbrace{
v_{\Lambda_1} \otimes \cdots \otimes v_{\Lambda_1}}_{(k-i) \; \;
\mbox{times}})= a x_{\alpha_1}(-1)^i \cdot v_{i \Lambda_0+(k-i)
\Lambda_1}, \nonumber
\end{eqnarray}
where $a$ a nonzero constant (note that $x_{\alpha_1}(-1) \cdot
v_{\Lambda_1}=0$ and $x_{\alpha_1}(-1)^2 \cdot v_{\Lambda_0}=0$).
Thus by (\ref{product oper}) and (\ref{hwc}) it follows that
\begin{eqnarray} \label{formula1}
&& e^{\lambda}_{(k)} \left ( x_{\alpha_1}(m_{1,1}) \cdots
x_{\alpha_1}(m_{r_1, 1}) \cdot v_{i \Lambda_1+(k-i) \Lambda_2} \right
) \\ &&= a x_{\alpha_1}(m_{1,1}-1) \cdots x_{\alpha_1}(m_{r_1,1}-1)
x_{\alpha_1}(-1)^i \cdot v_{i \Lambda_0+(k-i) \Lambda_1} \nonumber
\end{eqnarray} 
and
\begin{eqnarray} \label{formula2}
&& e^{\lambda}_{(k)} \left ( x_{\alpha_2}(m_{1,2}) \cdots
x_{\alpha_2}(m_{r_2, 2}) \cdot v_{i \Lambda_1+(k-i) \Lambda_2} \right
) \\ &&= a x_{\alpha_2}(m_{1,1}+1) \cdots x_{\alpha_2}(m_{r_2,2}+1)
x_{\alpha_1}(-1)^i \cdot v_{i \Lambda_0+(k-i) \Lambda_1}, \nonumber
\end{eqnarray} 
where $a \neq 0$, $r_j >0$, $m_{1,j}, \dots, m_{r_j, j} \in
\mathbb{Z}$ and $j=1,2$.

For $\beta = \alpha_2-\lambda_2=-\frac{1}{3} \alpha_1+\frac{1}{3}
\alpha_2$ we have the isomorphism
$$ e^{\beta}_{(k)}: V_P^{\otimes k} \longrightarrow V_P^{\otimes k}.
$$
As before, we obtain a linear map between principal subspaces
\begin{equation} \label{map1'}
e^{\beta}_{(k)}: W((k-i) \Lambda_1+i \Lambda_2) \longrightarrow W(i
\Lambda_0+(k-i) \Lambda_2)
\end{equation}
for any $i$ with $0 \leq i \leq k$. Thus it follows that
\begin{eqnarray} \label{formula'}
&&e^{\beta}_{(k)} ( x_{\alpha_j}(m_{1,j}) \cdots x_{\alpha_j}(m_{r_j,
j}) \cdot v_{(k-i)\Lambda_1+i \Lambda_2}) \\ && = a
x_{\alpha_j}(m_{1,j}-\langle \beta, \alpha_j \rangle) \cdots
x_{\alpha_j} (m_{r_j, j}-\langle \beta, \alpha_j \rangle) x_{\alpha_2}
(-1)^i \cdot v_{i \Lambda_0+(k-i) \Lambda_2} \nonumber
\end{eqnarray} 
for $a \neq 0$ (cf. (\ref{hwc}), (\ref{formula1}) and
(\ref{formula2})).

When $k=i=1$ the maps $e^{\lambda}_{(1)}$ and $e^{\beta}_{(1)}$ are
exactly the maps $e^{\lambda^{1}}$ and $e^{\lambda^{2}}$ defined in
Section 4 of \cite{C1}.

Intertwining operators among standard
$\widehat{{\goth{sl}(3)}}$-modules and fusion rules (the dimensions of
the vector spaces of intertwining operators of a certain type) are
important tools in our work.  In this paper we follow \cite{DL} for
the construction of distinguished intertwining operators and for some
of their relevant properties.

Let $i$ be an integer with $1 \leq i \leq k$. By using Chapter 13 of
\cite{DL} there exists a nonzero intertwining operator
 \begin{equation}
 {\cal Y} (v_{(k-1)\Lambda_0+\Lambda_1}, x):
 L(i\Lambda_0+(k-i)\Lambda_1) \longrightarrow
 L((i-1)\Lambda_0+(k-i+1)\Lambda_1)\{ x\}
 \end{equation} 
 of type
\begin{equation} \label{inter}
\left(
\begin{array}{c}
L((i-1) \Lambda_0 + (k-i+1) \Lambda_1)                           \\
\begin{array}{cc}
L((k-1)\Lambda_0+ \Lambda_1)   &  L( i \Lambda_0 + (k-i)\Lambda_1)
\end{array} 
\end{array}
\right ).
\end{equation}  
In fact, the dimension of the vector space of the intertwining
operators of type (\ref{inter}) (i.e. fusion rule) is one
(cf. \cite{DL}, \cite{FZ}, \cite{Li3}; see also \cite{BMW}, \cite{F}
and \cite{FW}). In this work we use only intertwining operators whose
fusion rules equal one.
   
We use the notation $$h_1=h_{(k-1)\Lambda_0+\Lambda_1}, \; \;
h_2=h_{i\Lambda_0+ (k-i)\Lambda_1}\; \; \mbox{and} \; \;
h_3=h_{(i-1)\Lambda_0+(k-i+1)\Lambda_1}
$$ 
(recall (\ref{h})). We have
\[
{\cal Y}(v_{(k-1) \Lambda_0+ \Lambda_1}, x) \in x^{h_3-h_1-h_2}
(\mbox{Hom} \; (L(i\Lambda_0+(k-i)\Lambda_1),
L((i-1)\Lambda_0+(k-i+1)\Lambda_1))) [[x, x^{-1}]]
\]
(cf. \cite{FHL}, \cite{FZ} and \cite{Li3}).
As a consequence of the Jacobi identity we obtain
\begin{equation}\label{brac}
[{\cal Y}(v_{(k-1) \Lambda_0+\Lambda_1}, x), U(\bar{\goth{n}})]=0,
\end{equation}
so that the coefficients of ${\cal Y}(v_{(k-1) \Lambda_0+\Lambda_1},
x)$ commute with the action of $U(\bar{\goth{n}})$ (cf. \cite{DL}; see
also \cite{G1}, \cite{CLM1}-\cite{CLM2} and \cite{C1}).
 
Now we take the constant term of $x^{-h_3+h_1+h_2} {\cal
Y}(v_{(k-1)\Lambda_0+\Lambda_1}, x)$ and we denote it by ${\cal
Y}_c(v_{(k-1)\Lambda_0+\Lambda_1}, x)$. The restriction of ${\cal
Y}_c(v_{(k-1)\Lambda_0+\Lambda_1}, x)$ to $W(i \Lambda_0+(k-i)
\Lambda_1)$ is a linear map
\begin{equation} \label{map2}
{\cal Y}_c(v_{(k-1)\Lambda_0+\Lambda_1}, x): W(i \Lambda_0+(k-i)
\Lambda_1) \longrightarrow W((i-1)\Lambda_0+ (k-i+1) \Lambda_1)
\end{equation} 
with the following property:
\begin{equation} \label{p1}
{\cal Y}_c(v_{(k-1)\Lambda_0+\Lambda_1}, x) \; (v_{i \Lambda_0+(k-i)
\Lambda_1})= bv_{(i-1)\Lambda_0+ (k-i+1) \Lambda_1},
\end{equation} 
where $b$ is a nonzero scalar. By (\ref{brac}) we immediately have
\begin{equation} \label{p2}
[{\cal Y}_c(v_{(k-1)\Lambda_0+\Lambda_1}, x), U(\bar{\goth{n}})]= 0,
\end{equation}
so that ${\cal Y}_c(v_{(k-1)\Lambda_0+\Lambda_1}, x)$ commutes with
the action of $U(\bar{\goth{n}})$.

Similarly, there exists a nonzero intertwining operator of type 
\begin{equation} \label{inter1}
\left(
\begin{array}{c}
L((i-1) \Lambda_0 + (k-i+1) \Lambda_2)                           \\
\begin{array}{cc}
L((k-1)\Lambda_0+ \Lambda_2)   &  L( i \Lambda_0 + (k-i)\Lambda_2)
\end{array} 
\end{array}
\right ),
\end{equation}
and there also exists a linear map associated to the intertwining
operator ${\cal Y}(v_{(k-1)\Lambda_0+\Lambda_2}, x)$ of type
(\ref{inter1}),
\begin{equation} \label{map2'}
{\cal Y}_c(v_{(k-1)\Lambda_0+\Lambda_2}, x):
W(i\Lambda_0+(k-i)\Lambda_2) \longrightarrow
W((i-1)\Lambda_0+(k-i+1)\Lambda_2),
\end{equation}
satisfying similar properties to those of the map (\ref{map2}).

The standard modules $L(k \Lambda_0)$, $L(k\Lambda_1)$ and
$L(k\Lambda_2)$ are ``group-like'' elements in the fusion ring at
level $k$. Such modules are also called simple currents
(cf. \cite{Li2}).

When $k=i=1$, the maps (\ref{map2}) and (\ref{map2'}) are the maps
${\cal Y}_c(e^{\lambda_1},x)$ and ${\cal Y}_c(e^{\lambda_2}, x)$ used
in \cite{C1}.

We now prove our main theorem that gives two families of $i$ exact
sequences of maps between principal subspaces for $1 \leq i \leq k$.

\begin{theorem} \label{theorem1}
Consider the maps $e^{\lambda}_{(k)}$, $e^{\beta}_{(k)}$, ${\cal
Y}_c(v_{(k-1)\Lambda_0+\Lambda_1}, x)$ and ${\cal Y}_c
(v_{(k-1)\Lambda_0+\Lambda_2}, x)$ introduced above (recall
(\ref{map1}), (\ref{map1'}), (\ref{map2}) and (\ref{map2'})).  Then
for any $i$ with $1 \leq i \leq k$ the following sequences:
\begin{eqnarray} \label{sequence}
\lefteqn{0 \longrightarrow W(i \Lambda_1 + (k-i) \Lambda_2)
\stackrel{e^{\lambda}_{(k)}} \longrightarrow} \\ && W(i \Lambda_0 +
(k-i) \Lambda_1) \stackrel{{\cal Y}_{c}(v_{(k-1)\Lambda_0+\Lambda_1},
x)} \longrightarrow \nonumber \\ && \hspace{2em} W((i-1)
\Lambda_0+(k-i+1) \Lambda_1) \longrightarrow 0 \nonumber
\end{eqnarray}
and
\begin{eqnarray} \label{sequence'}
\lefteqn{0 \longrightarrow W((k-i) \Lambda_1 + i \Lambda_2)
\stackrel{e^{\beta}_{(k)}} \longrightarrow} \\ && W(i \Lambda_0 +
(k-i) \Lambda_2) \stackrel{{\cal Y}_{c}(v_{(k-1)\Lambda_0+\Lambda_2},
x)} \longrightarrow \nonumber \\ && \hspace{2em} W((i-1)
\Lambda_0+(k-i+1) \Lambda_2) \longrightarrow 0 \nonumber
\end{eqnarray}
are exact.
\end{theorem}
{\em Proof:} We will prove that the sequence (\ref{sequence}) is
exact. The proof of the exactness of the sequence (\ref{sequence'}) is
completely analogous and we omit it.

We already know by (\ref{formula1}) and (\ref{formula2}) that
$e^{\lambda}_{(k)}$ maps $W(i \Lambda_1+(k-i) \Lambda_2)$ to $W(i
\Lambda_0+(k-i) \Lambda_1)$. This map is clearly injective. Since
${\cal Y}_c(v_{(k-1)\Lambda_0+\Lambda_1}, x)$ maps the highest weight
vector $v_{i\Lambda_0+(k-i)\Lambda_1}$ to a nonzero multiple of the
highest weight vector $v_{(i-1)\Lambda_0+(k-i+1) \Lambda_1}$ (recall
(\ref{p1})) and since this map commutes with the action of
$U(\bar{\goth{n}})$ (recall (\ref{p2})) we have that ${\cal
Y}_c(v_{(k-1)\Lambda_0+\Lambda_1}, x)$ is surjective.

Let $w \in \mbox{Im} \; e^{\lambda}_{(k)}$. By (\ref{formula1}) and
(\ref{formula2}) we have
$$
w=vx_{\alpha_1}(-1)^i \cdot v_{i \Lambda_0+(k-i) \Lambda_1}
$$
with $v \in U(\bar{\goth{n}})$. Hence
\begin{equation} \label{a}
{\cal Y}_c(v_{(k-1)\Lambda_0+\Lambda_1}, x) (w)= b v
x_{\alpha_1}(-1)^i \cdot v_{(i-1) \Lambda_0+(k-i) \Lambda_1},
\end{equation}
where $b$ is a nonzero constant (cf. (\ref{p1}) and (\ref{p2})). Now
(\ref{a}) combined with
$$ x_{\alpha}(-1) \cdot v_{\Lambda_1}=0 \; \; \mbox{and} \; \;
x_{\alpha}(-1)^2 \cdot v_{\Lambda_0}=0
$$
implies
$$
{\cal Y}_c(v_{(k-1)\Lambda_0+\Lambda_1}, x) (w)=0,
$$
and thus it gives the inclusion
\begin{equation} \label{i-k}
\mbox{Im} \; e^{\lambda}_{(k)} \subset \mbox{Ker} \; {\cal
Y}_c(v_{(k-1) \Lambda_0+\Lambda_1}, x).
\end{equation}

It remains to prove the inclusion
\begin{equation}
\mbox{Ker} \; {\cal Y}_c(v_{(k-1) \Lambda_0+\Lambda_1}, x) \subset
\mbox{Im} \; e^{\lambda}_{(k)}.
\end{equation}
In order to show this inclusion we first charaterize the vector spaces
$\mbox{Ker} \; {\cal Y}_c(v_{(k-1)\Lambda_0+\Lambda_1}, x)$ and
$\mbox{Im}\; e^{\lambda}_{(k)}$. Let $w \in \mbox{Ker} \; {\cal
Y}_c(v_{(k-1)\Lambda_0+\Lambda_1}, x)$. In particular $w \in
W(i\Lambda_0+(k-i)\Lambda_1)$, so that 
$$
w =f_{i \Lambda_0+(k-i) \Lambda_1} (u)
$$ for $u \in U(\bar{\goth{n}})$ (cf. (\ref{map})). By using again
(\ref{p1}) and (\ref{p2}) we obtain
$$ {\cal Y}_c(v_{(k-1)\Lambda_0+\Lambda_1}, x) (f_{i \Lambda_0+(k-i)
\Lambda_1} (u))=0 \Leftrightarrow f_{(i-1)\Lambda_0+(k-i+1)\Lambda_1}
(u)=0.
$$
We have just shown that
\begin{equation} \label{chara1}
w=f_{i \Lambda_0+(k-i) \Lambda_1} (u) \in \mbox{Ker} \; {\cal
Y}_c(v_{(k-1)\Lambda_0+\Lambda_1}, x) \Leftrightarrow u \in
I_{(i-1)\Lambda_0+(k-i+1) \Lambda_1}.
\end{equation}

Let $w \in \mbox{Im} \; e^{\lambda}_{(k)}$. Then by (\ref{formula1})
and (\ref{formula2}) we have
\begin{equation} \label{e1}
w=vx_{\alpha_1}(-1)^i \cdot v_{i\Lambda_0+(k-i)\Lambda_1} = f_{i
\Lambda_0+(k-i)\Lambda_1}(v x_{\alpha_1}(-1)^i),
\end{equation}
where $v \in U(\bar{\goth{n}})$. On the other hand, since $w \in
\mbox{Im} \; e^{\lambda}_{(k)} \subset W(i\Lambda_0+(k-i) \Lambda_1)$,
\begin{equation} \label{e2}
w=f_{i \Lambda_0+(k-i)\Lambda_1}(u)
\end{equation}
with $u \in U(\bar{\goth{n}})$.
Now (\ref{e1}) and (\ref{e2}) imply
$$
u-vx_{\alpha_1}(-1)^i \in I_{i \Lambda_0+(k-i)\Lambda_1}.
$$
Thus we have obtained
\begin{equation} \label{chara2}
w=f_{i \Lambda_0+(k-i)\Lambda_1}(u) \in \mbox{Im} \; e^{\lambda}_{(k)}
\Leftrightarrow u \in I_{i\Lambda_0+(k-i) \Lambda_1} +
U(\bar{\goth{n}})x_{\alpha_1}(-1)^i.
\end{equation}

We also have
\begin{equation} \label{chara3}
I_{(i-1)\Lambda_0+(k-i+1) \Lambda_1}= I_{k \Lambda_0}
+U(\bar{\goth{n}}) x_{\alpha_1}(-1)^i
\end{equation}
and
\begin{equation}\label{chara4}
I_{i\Lambda_0+(k-i) \Lambda_1}= I_{k \Lambda_0} + U(\bar{\goth{n}})
x_{\alpha_1}(-1)^{i+1},
\end{equation}
so that
\begin{equation} \label{chara5}
I_{(i-1)\Lambda_0+(k-i+1)\Lambda_1}= I_{k\Lambda_0}
+U(\bar{\goth{n}}_{-}) x_{\alpha_1}(-1)^i \subset
I_{i\Lambda_0+(k-i)\Lambda_1} + U(\bar{\goth{n}}_{-})
x_{\alpha_1}(-1)^i
\end{equation} 
(recall Corollary \ref{discrepancy}).  Now (\ref{chara1}),
(\ref{chara2}) and (\ref{chara5}) prove the inclusion
\begin{equation}
\mbox{Ker} \; {\cal Y}_c(v_{(k-1)\Lambda_0+\Lambda_1}, x) \subset
\mbox{Im} \; e^{\lambda}_{(k)},
\end{equation}
and thus the exactness of the sequence (\ref{sequence}). This
concludes our theorem. $\; \; \; \Box$
\vspace{1em}

\begin{remark}
\rm The chain property of the sequence (\ref{sequence}) follows from
(\ref{chara1}), (\ref{chara2}), (\ref{chara3}) and (\ref{chara4}). The
chain property follows also by a different argument, as we have just
seen in the proof of the previous theorem. We do not have an
elementary proof that does not use the description of the ideals
$I_{i\Lambda_0+(k-i)\Lambda_1}$ and
$I_{(i-1)\Lambda_0+(k-i+1)\Lambda_1}$ for the exactness of the
sequence (\ref{sequence}). Similar assertions are of course true for
the sequence (\ref{sequence'}).
\end{remark}

Recall from Lemmas \ref{lemma1} and \ref{lemma2} the equations
\begin{equation} \label{1-0-new}
\chi'_{W(k \Lambda_1)}(x_1, x_2; q)= \chi'_{W(k \Lambda_0)}(x_1q, x_2;
q)
\end{equation}
and
\begin{equation} \label{2-0-new}
\chi'_{W(k \Lambda_2)}(x_1, x_2; q)=\chi'_{W(k\Lambda_0)}(x_1, x_2q;
q).
\end{equation}

Given the exact sequences from Theorem
\ref{theorem1} we now derive a system of $q$-difference equations that
characterize the graded dimensions of $W(i\Lambda_0+(k-i)\Lambda_j)$
for $1 \leq i \leq k$ and $j=1,2$.

\begin{theorem} \label{q-equations}
For any $i$ such that $1 \leq i < k$ we have the following
$q$-difference equations:
\begin{eqnarray} \label{main rec}
&& \chi'_{W(i\Lambda_0+(k-i)\Lambda_1)}(x_1, x_2;q)-
\chi'_{W((i-1)\Lambda_0+(k-i+1)\Lambda_1)}(x_1, x_2; q) \\ \nonumber
&& +x_1^ix_2^{i-k}q^k \chi'_{W((k-i-1)\Lambda_0+(i+1)\Lambda_2)}
(x_1q^2, x_2q^{-2}; q) - x_1^ix_2^{i-k}q^k
\chi'_{W((k-i)\Lambda_0+i\Lambda_2)}(x_1q^2, x_2q^{-2}; q) \\
\nonumber && \hspace{4em} =0 \nonumber
\end{eqnarray}
and
\begin{eqnarray} \label{main recur}
&& \chi'_{W(i\Lambda_0+(k-i)\Lambda_2)}(x_1, x_2;q)-
\chi'_{W((i-1)\Lambda_0+(k-i+1)\Lambda_2)}(x_1, x_2; q) \\ \nonumber
&& +x_1^{i-k}x_2^{i}q^k \chi'_{W((k-i-1)\Lambda_0+(i+1)\Lambda_1)}
(x_1q^{-2}, x_2q^{2}; q) - x_1^{i-k}x_2^{i}q^k
\chi'_{W((k-i)\Lambda_0+i\Lambda_1)}(x_1q^{-2}, x_2q^{2}; q) \\
\nonumber && \hspace{4em} =0. \nonumber
\end{eqnarray}
If $i=k \geq 1$ then we have
\begin{eqnarray} \label{rec1}
&& \chi'_{W((k-1)\Lambda_0+\Lambda_1)}(x_1, x_2; q)
-\chi'_{W((k-1)\Lambda_0+\Lambda_2)}(x_1, x_2; q) \\ \nonumber && +
(x_1q)^k\chi'_{W(k \Lambda_1)}(x_1q, x_2q^{-1}; q) -(x_2q)^k\chi'_{W(k
\Lambda_2)}(x_1q^{-1}, x_2q; q)=0 .\nonumber
\end{eqnarray}
\end{theorem}
{\em Proof:} We first prove that the following $q$-difference
equations hold:
\begin{eqnarray} \label{recursion1}
&&\chi'_{W(i\Lambda_0+(k-i)\Lambda_1)}(x_1, x_2;q) \\ \nonumber
&&=x_1^iq^i \chi'_{W(i\Lambda_1+(k-i)\Lambda_2)}(x_1q, x_2q^{-1}; q)+
\chi'_{W((i-1)\Lambda_0+(k-i+1)\Lambda_1)}(x_1, x_2; q) \nonumber
\end{eqnarray}
and
\begin{eqnarray} \label{recursion2}
&& \chi'_{W(i\Lambda_0+(k-i)\Lambda_2)} (x_1, x_2; q)\\ \nonumber && =
x_2^iq^i \chi'_{W((k-i)\Lambda_1+i\Lambda_2)}(x_1q^{-1}, x_2q;
q)+\chi'_{W((i-1)\Lambda_0+(k-i+1)\Lambda_2)}(x_1, x_2; q). \nonumber
\end{eqnarray}
The $q$-equation (\ref{recursion1}) follows by using the exact
sequence (\ref{sequence}). Formula (\ref{recursion2}) is obtained
similarly from the exactness of the sequence (\ref{sequence'}).

Let us consider the homogeneous subspace $W(i \Lambda_1+(k-i)
\Lambda_2)'_{r_1, r_2; s}$ (recall our notation (\ref{w'})). From
(\ref{formula1}) and (\ref{formula2}) we notice that the map
$e^{\lambda}_{(k)}$ increases the charge corresponding to $\lambda_1$
by $i$ and preserves the charge corresponding to $\lambda_2$. Also,
for any $w \in W(i\Lambda_1+(k-i)\Lambda_2)'_{r_1, r_2; s}$ the vector
$e^{\lambda}_{(k)}(w)$ has weight $s+r_1-r_2+i$.  Hence we have
\begin{equation} 
e^{ \lambda}_{(k)}: W(i \Lambda_1 + (k-i) \Lambda_2)'_{ r_1, r_2;s}
\longrightarrow W(i \Lambda_0 + (k-i) \Lambda_1)'_{r_1+i, r_2;
s+r_1-r_2+i}
\end{equation}
or equivalently,
\begin{equation} \label{rel}
e^{ \lambda}_{(k)}: W(i \Lambda_1 + (k-i) \Lambda_2)'_{ r_1-i, r_2;
s-r_1+r_2} \longrightarrow W(i \Lambda_0 + (k-i) \Lambda_1)'_{r_1,
r_2; s}.
\end{equation}

By (\ref{p1}) and (\ref{p2}) we obtain 
\begin{equation}
{\cal Y}_c(v_{(k-1)\Lambda_0+ \Lambda_1}, x): W(i \Lambda_0 + (k-i)
\Lambda_1)'_{r_1, r_2; s} \longrightarrow W((i-1)\Lambda_0 +
(k-i+1)\Lambda_1)'_{r_1, r_2; s}.
\end{equation}
The exacteness of the sequence (\ref{sequence}) implies the equation
(\ref{recursion1}) for any $i$ with $1 \leq i \leq k$.

Easy computations show that (\ref{recursion1}) and (\ref{recursion2})
yield the following equations:
\begin{eqnarray} \label{rec1-new}
&&\chi'_{W(i\Lambda_1+(k-i)\Lambda_2)}(x_1, x_2; q) \\ \nonumber && =
x_1^{-i} \chi'_{W(i\Lambda_0+(k-i)\Lambda_1)}(x_1q^{-1}, x_2q; q)-
x_1^{-i} \chi'_{W((i-1)\Lambda_0+(k-i+1)\Lambda_1)}(x_1q^{-1}, x_2q;
q) \nonumber
\end{eqnarray}
and
\begin{eqnarray} \label{rec2-new}
&& \chi'_{W((k-i)\Lambda_1+i\Lambda_2)}(x_1, x_2; q) \\ \nonumber && =
x_2^{-i} \chi'_{W(i\Lambda_0+(k-i)\Lambda_2)} (x_1q, x_2q^{-1};
q)-x_2^{-i}\chi'_{W((i-1)\Lambda_0+(k-i+1)\Lambda_2)}(x_1q, x_2q^{-1};
q) \nonumber
\end{eqnarray}
for all $i$ with $1 \leq i \leq k$. 

Now assume that $1 \leq i <k$. We take $i$ instead of $k-i$ in
(\ref{rec2-new}) and equate the expression of
$\chi'_{W(i\Lambda_1+(k-i)\Lambda_2)}(x_1, x_2;q)$ from the resulting
formula with the expression of
$\chi'_{W(i\Lambda_1+(k-i)\Lambda_2)}(x_1, x_2; q)$ from
(\ref{rec1-new}) and thus we obtain
\begin{eqnarray} \nonumber
&& x_1^{-i} \chi'_{W(i\Lambda_0+(k-i)\Lambda_1)}(x_1q^{-1}, x_2q;
q)-x_1^{-i}\chi'_{W((i-1)\Lambda_0+(k-i+1)\Lambda_1)}(x_1q^{-1}, x_2q;
q) \nonumber \\ &&= x_2^{-k+i}
\chi'_{W((k-i)\Lambda_0+i\Lambda_2)}(x_1q, x_2q^{-1};
q)-x_2^{-k+i}\chi'_{W(k-i-1)\Lambda_0+(i+1)\Lambda_2)}(x_1q,
x_2q^{-1}; q). \nonumber
\end{eqnarray}
By multiplying all terms of this expression by $x_1^i$ and substituing
 $x_1$ with $x_1q$ and $x_2$ with $x_2q^{-1}$ we obtain the four-term
 $q$-difference equation (\ref{main rec}).

In a similar way, we obtain (\ref{main recur}) by taking $k-i$ instead
of $i$ in (\ref{rec1-new}), equating the expressions of
$\chi'_{W((k-i)\Lambda_1+i\Lambda_2)}(x_1, x_2; q)$ and doing the
corresponding substitutions for $x_1$ and $x_2$. Note that (\ref{main
recur}) can be derived from (\ref{main rec}) by using an isomorphism
induced from the Dynkin diagram automorphism of $\goth{sl}(3)$ that
interchanges $\alpha_1$ and $\alpha_2$, and thus it interchanges
$\Lambda_1$ and $\Lambda_2$.

We now assume that $k=i$ in (\ref{rec1-new}) and (\ref{rec2-new}). We
first multiply (\ref{rec1-new}) by $x_1^k$ and (\ref{rec2-new}) by
$x_2^k$ and then we make the substitutions $x_1\mapsto x_1q$, $x_2
\mapsto x_2q^{-1}$ in (\ref{rec1-new}) and $x_1\mapsto x_1q^{-1}$,
$x_2 \mapsto x_2q$ in (\ref{rec2-new}). By eliminating the term
$\chi'_{W(k\Lambda_0)}(x_1, x_2; q)$ we obtain (\ref{rec1}).  $\; \;
\; \Box$
\vspace{1em}

The particular case $k=i=1$ of Theorem \ref{q-equations} is listed
below.
 
\begin{remark}
\rm By taking $k=i=1$ in (\ref{recursion1}) and (\ref{recursion2}) and by
using (\ref{1-0-new}) and (\ref{2-0-new}) for $k=1$ we recover the
recursions satisfied by the graded dimension of the principal subspace
$W(\Lambda_0)$ of the level 1 vacuum representation of 
$\widehat{\goth{sl}(3)}$ obtained in \cite{C1} (see also \cite{C2}):
\begin{equation}
\chi_{W(\Lambda_0)}(x_1,x_2;q)=\chi_{W(\Lambda_0)}(x_1q,x_2;q)
+x_1q\chi_{W(\Lambda_0)}(x_1q^2, x_2q^{-1};q)
\end{equation}
and
\begin{equation}
\chi_{W(\Lambda_0)}(x_1,x_2;q)=\chi_{W(\Lambda_0)}(x_1,x_2q;q)
+x_2q\chi_{W(\Lambda_0)}(x_1q^{-1}, x_2q^{2};q)
\end{equation}
(recall from \cite{C1} that $\chi'_{W(\Lambda_0)}(x_1, x_2;
q)=\chi_{W(\Lambda_0)}(x_1, x_2; q)$). These recursions are equivalent
with
\begin{eqnarray} 
&& \chi_{W(\Lambda_0)}(x_1q, x_2; q)-\chi_{W(\Lambda_0)}(x_1, x_2q; q)
\nonumber \\ && +x_1q\chi_{W(\Lambda_0)}(x_1q^2, x_2q^{-1}; q)- x_2q
\chi_{W(\Lambda_0)}(x_1q^{-1}, x_2q^2; q)=0. \nonumber
\end{eqnarray}
\end{remark}

Graded dimensions of the principal subspaces
 $W(i\Lambda_0+(k-i)\Lambda_j)$ with $0 \leq i \leq k$ and $j=1,2$
 were computed in \cite{G1} by constructing quasiparticle bases of
 these principal subspace:
 \begin{eqnarray} \label{grd}
&& \chi'_{W(i\Lambda_0+(k-i) \Lambda_j)}(x_1, x_2;q) \\ \nonumber && =
 \sum_{{\begin{array}{c} 0 \leq M_k \leq \cdots \leq M_1, \\ 0 \leq
 N_k \leq \cdots \leq N_1 \end{array}}} \frac{q^{\sum_{t=1}^k (M_t^2 +
 N_t^2- M_tN_t) +\sum_{t=1}^k (M_t \delta_{1, j_{t}} + N_t \delta_{2,
 j_{t}})}} {(q)_{M_1-M_2} \cdots (q)_{M_k} (q)_{N_1-N_2} \cdots
 (q)_{N_k}} x_1^{\sum_{t=1}^k M_t} x_2^{\sum_{t=1}^k N_t}, \nonumber
 \end{eqnarray} 
where $j_t=0$ for $0 \leq t \leq i$ and $j_t=j$ for $i<t \leq k$,
$j=1, 2$, and where we use the notation
$$
(q)_m=(1-q) \cdots (1-q^m)
$$ for any nonnegative integer $m$.  Following an argument similar to
the one developed in Section 7.3 in \cite{A} one can show that
(\ref{grd}) satisfies the system of equations (\ref{1-0-new}),
(\ref{2-0-new}), (\ref{main rec}), (\ref{main recur}) and
(\ref{rec1}).

The $q$-difference equations (\ref{rec1-new}) and (\ref{rec2-new})
combined with formulas (\ref{grd}) give the graded dimensions of the
principal subspaces corresponding to highest weights of the form $i
\Lambda_1+(k-i)\Lambda_2$ for any $i$ with $1 \leq i \leq k$.

We first observe that by taking $k=i$ in (\ref{rec1-new}) and using
the formulas for $\chi'_{W(k\Lambda_0)}(x_1q^{-1}, x_2q; q)$ and
$\chi'_{W((k-1)\Lambda_0+\Lambda_1)}(x_1q^{-1}, x_2q; q)$ we recover
the graded dimension $\chi'_{W(k\Lambda_1)}(x_1, x_2; q)$.

Now assume that $1 \leq i \leq k-1$ in (\ref{rec1-new}).

\begin{corollary} \label{new} 
We have
\begin{eqnarray} \label{new grd}
&&\hspace{4em} \chi'_{W(i\Lambda_1+(k-i)\Lambda_2)} (x_1, x_2;q)= \\
&& \sum_{{\begin{array}{c} 0 \leq M_k\leq \cdots \leq M_1, \\
0\leq N_k \leq \cdots \leq N_1 \end{array}}} \frac{ q^{ \sum_{t=1}^k
(M_t^2+N_t^2-M_tN_t) + \sum_{t=i+1}^k M_t }(1-q^{M_i}) q^{\sum_{t=1}^k
(N_t-M_t)} }{(q)_{M_1-M_2} \cdots (q)_{M_k} (q)_{N_1-N_2} \cdots
(q)_{N_k}} \nonumber \\ && \hspace{10em}
\sum_{{\begin{array}{c} 0 \leq M_k\leq \cdots \leq M_1, \\ 0\leq
N_k \leq \cdots \leq N_1
\end{array}}}x_1^{\sum_{t=1}^k M_t}x_1^{-i} x_2^{\sum_{t=1}^k
N_t}. \nonumber
\end{eqnarray}
\end{corollary}
{\em Proof:} The statement follows by using (\ref{rec1-new}) together
with formula (\ref{grd}) for the graded dimensions of
$W(i\Lambda_0+(k-i)\Lambda_1)$ and
$W((i-1)\Lambda_0+(k-i+1)\Lambda_1)$. $\; \; \; \Box$
\vspace{1em}

It is interesting to consider the problem of computing the graded
dimension of $W(k_0\Lambda_0+k_1\Lambda_1+k_2\Lambda_2)$, where $k_0,
k_1, k_2$ are any positive integers whose sum is $k$, as a combination
of the graded dimensions already obtained.  One attempt is given by
the next theorem, whose proof is completely analogous to that of
Theorem \ref{theorem1}.
 
Let $k_0, k_1, k_2 \in \mathbb{N}$ such that $k_1, k_2 \geq 1$ and
$k_0+k_1+k_2=k$. We consider intertwining operators of types
\begin{equation} \label{intert-1-new}
\left(
\begin{array}{c}
L((k_1-1) \Lambda_0 + (k_0+k_2+1) \Lambda_1)                           \\
\begin{array}{cc}
L((k_1+k_2-1)\Lambda_0+ (k_0+1)\Lambda_1) & L( (k_0+k_1) \Lambda_0 +
k_2\Lambda_1)
\end{array} 
\end{array}
\right )
\end{equation} 
and
\begin{equation} \label{intert-2-new}
\left(
\begin{array}{c}
L((k_2-1) \Lambda_0 + (k_0+k_1+1) \Lambda_2)                           \\
\begin{array}{cc}
L((k_1+k_2-1)\Lambda_0+ (k_0+1)\Lambda_2) & L( (k_0+k_2) \Lambda_0 +
k_1\Lambda_2)
\end{array} 
\end{array}
\right ).
\end{equation} 
By \cite{DL} (see also \cite{FZ} and \cite{Li3}) we have that the
dimension of the vector spaces of intertwining operators of type
(\ref{intert-1-new}) and of type (\ref{intert-2-new}), respectively,
equals one. Consider the constant terms of the intertwining operators
${\cal Y}(v, x)$ of type (\ref{intert-1-new}) and ${\cal Y}(v', x)$ of
type (\ref{intert-2-new}):
$$ {\cal Y}_c(v, x): W((k_0+k_1)\Lambda_0+k_2\Lambda_1)
\longrightarrow W((k_1-1)\Lambda_0+(k_0+k_2+1)\Lambda_1)
$$
and
$$ {\cal Y}_c (v', x): W((k_0+k_2)\Lambda_0+k_1\Lambda_2)
\longrightarrow W((k_2-1)\Lambda_0+(k_0+k_1+1)\Lambda_2),
$$ where $v=v_{(k_1+k_2-1)\Lambda_0+(k_0+1)\Lambda_1}$ and
$v'=v_{(k_1+k_2-1)\Lambda_0+(k_0+1)\Lambda_2}$ are highest weight
vectors of $L((k_1+k_2-1) \Lambda_0+(k_0+1)\Lambda_1)$ and of
$L((k_1+k_2-1)\Lambda_0+(k_0+1)\Lambda_2)$.
 
Recall the linear maps (\ref{id}), (\ref{map1}) and
(\ref{map1'}). Consider
$$ \mbox{Id}_{(k_0)} \otimes e^{\lambda}_{(k-k_0)}:
W(k_0\Lambda_0+k_1\Lambda_1+k_2\Lambda_2) \longrightarrow
W((k_0+k_1)\Lambda_0+k_2\Lambda_1)
$$
and
$$ \mbox{Id}_{(k_0)} \otimes e^{\beta}_{(k-k_0)}:
W(k_0\Lambda_0+k_1\Lambda_1+k_2\Lambda_2) \longrightarrow
W((k_0+k_2)\Lambda_0+k_1\Lambda_2).
$$
\begin{theorem} \label{general}
There are natural exact sequences
\begin{eqnarray} 
\lefteqn{0 \longrightarrow W(k_0\Lambda_0+k_1\Lambda_1+k_2\Lambda_2)
\stackrel{\rm{Id}_{(k_0)} \otimes e^{\lambda}_{(k-k_0)}}
\longrightarrow} \\ && W((k_0+k_1) \Lambda_0+k_2\Lambda_1)
\stackrel{{\cal Y}_{c}(v, x)} \longrightarrow \nonumber \\ &&
\hspace{2em} W((k_1-1)\Lambda_0+(k_0+k_2+1)\Lambda_1) \longrightarrow
0 \nonumber
\end{eqnarray}
and
\begin{eqnarray} 
\lefteqn{0 \longrightarrow W(k_0\Lambda_0+k_1\Lambda_1+k_2\Lambda_2)
\stackrel{\rm{Id}_{(k_0)} \otimes e^{\beta}_{(k-k_0)}}
\longrightarrow} \\ && W((k_0+k_2) \Lambda_0+k_1\Lambda_2)
\stackrel{{\cal Y}_{c}(v', x)} \longrightarrow \nonumber \\ &&
\hspace{2em} W((k_2-1)\Lambda_0+(k_0+k_1+1)\Lambda_2) \longrightarrow
0.  \; \; \; \; \; \; \; \Box \nonumber
\end{eqnarray}

\end{theorem}

However, the maps $\mbox{Id}_{(k_0)} \otimes e^{\lambda}_{(k-k_0)}$
and $\mbox{Id}_{(k_0)} \otimes e^{\beta}_{(k-k_0)}$, restricted to the
corresponding homogeneous spaces, do not shift the weight of a
homogeneous element by a scalar, so that our method used before in
Theorem \ref{q-equations} to obtain $q$-difference equations does not
extend to this case.

\vspace{.4in}
  
\noindent {\small \sc Department of Mathematics, Rutgers University,
Piscataway, NJ 08854} 
\vspace{.1in}

\noindent Current address: 

\noindent{\small \sc Department of Mathematics,  Ohio State University,
Columbus, OH 43210} \\
{\em E--mail address}: calinescu@math.ohio-state.edu

 \end{document}